\documentclass{amsart}	
\usepackage{amsfonts,amsmath,amsthm,amssymb,graphicx,euscript,mathrsfs}
\usepackage[all]{xy}
\DeclareGraphicsExtensions{.eps}
\numberwithin{equation}{section}

\newtheorem{thrm}[equation]{Theorem}
\newtheorem*{mthrm}{Theorem 3.3}

\newcommand{\X}{\mathcal{X}}
\newcommand{\F}{\mathbb{F}}
\newcommand{\Fq}{\mathbb{F}_q}
\newcommand{\LL}{\mathcal{L}}
\newcommand{\vP}{\emph{v}_P}
\newcommand{\Z}{\mathbb{Z}}

\theoremstyle{definition}

\newtheorem{expl}[equation]{Example}

\newtheorem{rmrk}[equation]{Remark}

\newtheorem*{note}{Note}

\begin{document} 	

\title[Minimum Distance of Geometric Goppa Codes]  {A Generalized
Floor Bound for the Minimum Distance of Geometric Goppa Codes and its
application to Two-Point Codes} 
\author{Benjamin Lundell and Jason McCullough}
\date{\today}
\address{Dept.~of Mathematics, University of Illinois, Urbana, Illinois 61801}
\email{blundell@math.uiuc.edu, jmccullo@math.uiuc.edu}

\begin{abstract}

We prove a new bound for the minimum distance of geometric Goppa codes
that generalizes two previous improved bounds.  We include examples of
the bound to one and two point codes over both the Suzuki and
Hermitian curves.
\end{abstract}

\maketitle
\section{Introduction}

In \cite{Goppa}, Goppa gives a construction of a family of
error-correcting codes using two divisors on a curve.  He also gives a
distance bound based on the degrees of the divisors.  Yang et. al. in \cite{Yang} and Chen et. al. in \cite{DuursmaChen} compute the actual minimum distance of the one-point Hermitian and certain one-point Suzuki codes respectively, often showing that the actual distance is significantly higher than Goppa's designed distance.  In \cite{F-R}, Feng and Rao present a decoding algorithm for Goppa codes, which always decodes at least up to half the designed distance.  The new distance bound obtained from their algorithm is in general very good as few codes have minimum distance exceeding the Feng-Rao (F-R) distance.

Other efforts to
improve and generalize the distance bounds have had some success.  In \cite{Pellikaan},
Kirfel and Pelikaan generalize a result of Garcia et al. in
\cite{GKL} and improve the designed distance by taking advantage of
pairs of large gaps in the Weierstrass gap sequence at a point $P$.  In
\cite{Matthews2.10}, Maharaj et al. improve the minimum distance by
using the notion of the floor of a divisor to capitalize on one large
gap in a multi-point code.  We prove a new bound, which generalizes
the previous two.  It can use more than one gap in the gap sequence, but
also applies nicely to multi-point codes.  We state our main result below.

\begin{mthrm}[Asymmetric Floor Bound]
Let $\X$ be a curve over $K / \Fq$ of genus $g$.  Let $P_1,\dots,P_n$
be distinct rational points on $\X$.  Define $D := P_1 + \dots + P_n$.
Let $A$, $B$, and $Z$ be divisors with support outside of $D$ such that
$Z$ is effective, $\ell(A) = \ell(A - Z)$, and $\ell(B) = \ell(B + Z)$.
 Then, putting $G=A+B$ yields

\[
d(C_\Omega(D,G)) \geq \rm{deg}(G) - (2g-2) +\rm{deg}(Z).
\]
\end{mthrm}

For some codes, the asymmetric bound achieves parameters that were
unattainable by the previous bounds.  See Section 4 for specific
examples.  Section 2 gives definitions and notation that will be used
throughout the paper.  Section 3 gives a formal statement of all three
mentioned bounds, a proof of the asymmetric bound, and relationships
among the three.

\section{Preliminaries}
Unless otherwise noted, we follow the same definitions and notations
as in \cite{Matthews2.10}.  If $P$ is a rational point on a curve
$\X$ that is defined over $\Fq$ with function field $K$, then $\vP$ represents the discrete 
valuation
corresponding to $P$.  For a divisor $A$, we denote the support of $A$
as $\rm{Supp}(A)$.  If $B$ is another divisor, the we define the greatest
common divisor of $A$ and $B$ by
\[
\rm{gcd}(A,B):= \sum_{P} \min\{\vP(A),\vP(B)\}P.
\]
From $A$, we can create two vector spaces, one of rational functions from the
function field $K$, and one of rational differentials:
\[
\LL(A) := \{f \in \F : (f) + A \succcurlyeq 0 \} \cup \{0\}
\]
and
\[
\Omega(A) := \{\eta \in \Omega : (\eta) \succcurlyeq A\} \cup \{0\}.
\]
We denote the dimension of $\LL(A)$ over $\Fq$ by $\ell(A)$ and the
dimension of $\Omega(A)$ over $\Fq$ by $i(A)$. 
\\ \\
Using the distinct rational points $P_1, \dots, P_n,Q_1, \dots, Q_m$
on $\X$ to form the two divisors 
\[
D:=P_1+ \dots + P_n \text{ and } 
G:=
\sum_{i} \alpha_iQ_i \quad \alpha_i \in \Z,
\]
 we can define two linear $m$-point codes:
\[
C_{\LL}(D,G):=\{(f(P_1), \dots,f(P_n)):f \in \LL(G)\}
\]
and
\[
C_{\Omega}(D,G):=\{(\rm{res}_{P_1}(\eta), \dots,\rm{res}_{P_n}(\eta)):\eta \in 
\Omega(G-D)\}.
\]
Both of these codes have length $n$.  The dimension of $C_{\LL}$ is
$\ell(G)-\ell(G-D)$ and the designed distance is $n-\rm{deg}(G)$.  The
dimension of $C_{\Omega}$ is $i(G-D)-i(G)$ and the designed
distance is $\rm{deg}(G)-(2g-2)$.

The floor of $A$ is defined in \cite{Matthews2.10} as the
unique divisor, $A'$, of minimum degree such that $\LL(A')=\LL(A)$
and is denoted $\lfloor A \rfloor$.

Finally, we say that an integer $\alpha$ is an $A$-gap at a point $P$
if and only if $\LL(A+\alpha P)=\LL(A + (\alpha -1)P)$.

\section{Asymmetric Bound}
We start by stating two theorems.  The first appears as Theorem 2.10
in \cite{Matthews2.10}, the second as Proposition 3.10 in
\cite{Pellikaan}.
\begin{thrm}[Floor Bound] \label{floorbd}
Let $K/\Fq$ be a function field of genus g.  Let $D:=P_1+ \dots + P_n$
where $P_1, \dots, P_n$ are distinct rational places of $F$, and let
$G:=H + \lfloor H \rfloor$ be a divisor of $F$ such that the support
of $H$ does not contain any of the places $P_1, \dots, P_n$.  Set
$E_H:=H - \lfloor H \rfloor$.  Then $C_\Omega(D,G)$ is an $[n,k,d]$
code whose parameters satisfy
\[
d \geq \rm{deg}(G)-(2g-2)+\rm{deg}(E_H) =  2\rm{deg}(H)-(2g-2).
\]
\end{thrm}

\begin{thrm}[K-P Bound] \label{Prop3.10}
Suppose that each of the integers $\alpha, \alpha +1, \dots, \alpha +
t$ is an $F$-gap at $P$ and $\beta -t, \dots,\beta -1 , \beta$ are
$G$-gaps at P.  Put $H=F+G+(\alpha+\beta-1)P$.  Suppose $D:=P_1+ \dots
+ P_n$, where the $P_i$ are n distinct rational points, each not equal
to $P$ and not belonging to the support of $H$.  Then the minimum
distance of $C_\Omega(D,H)$ is at least $\rm{deg}(H)-(2g-2)+(t+1)$.

\end{thrm}

We adapt the proof of the floor bound to prove the following theorem, 
which generalizes both the floor bound and the K-P bound. 

\begin{thrm}[Asymmetric Floor Bound (AF Bound)] \label{asymbd}
Let $\X$ be a curve over $K / \Fq$ of genus $g$.  Let $P_1,\dots,P_n$
be distinct rational points on $\X$.  Define $D := P_1 + \dots + P_n$.
Let $A$, $B$, $G$, and $Z$ be divisors with support outside of $D$ such that
$Z$ is effective, $\ell(A) = \ell(A - Z)$, $\ell(B) = \ell(B + Z)$, and $G = A + B$.  Then

\[
d(C_\Omega(D,G)) \geq \rm{deg}(G) - (2g-2) +\rm{deg}(Z).
\]
\end{thrm}
 
\proof Let $\eta \in \Omega(G-D)$ such that the code word $\vec c :=
(\rm{res}_{P_1}(\eta), \dots,\rm{res}_{P_n}(\eta))$ is of minimum nonzero weight.
Without loss of generality, we may assume that $c_i \neq 0$ for $1
\leq i \leq d$ and $c_i = 0$ for $d < i \leq n$.  Let $D':= P_1 +
\dots P_d$.  Since $\vec c$ is zero outside $D'$, we must have that
$(\eta) \succcurlyeq G-D'$.  Thus, there exists an effective divisor,
$E$, with support disjoint from that of $D'$ so that $W := (\eta) = G
- D' + E$.  Taking degrees of both sides we get that 
\[
2g-2 = \rm{deg}(G) - d + \rm{deg}(E) \Rightarrow d = \rm{deg}(G) - (2g-2) + \rm{deg}(E).
\]
Observe that,
\[
\rm{deg}(E) \geq \ell(A + E) - \ell(A) = \ell(A + E) - \ell(A - Z) \geq \ell(A + E) 
- \ell(A + E - Z).
\]

By applying the Riemann-Roch Theorem twice, we get that

\[
\ell(A + E) = \rm{deg}(A+E) +1 -g + \ell(W-(A+E))
\]
and
\[
\ell(A + E - Z) = \rm{deg}(A+E-Z) +1-g+\ell(W-(A+E-Z))
\]
and thus,
\[
\ell(A + E)-\ell(A + E - Z) = \rm{deg}(Z)+\ell(B-D')-\ell(B+Z-D').
\]
Now 
\[ \LL(B+Z-D') \subseteq \LL(B+Z) = \LL(B) \] 
implies that 
\[ \LL(B + Z - D') = \LL(B+Z-D') \cap \LL(B) = \LL(\rm{gcd}(B+Z-D',B)).
\]

Since $B$ and $Z$ have support outside of $D'$ and since $Z$ is effective,  
\\ $\rm{gcd}(B+Z-D',B)=B-D'$, and thus $\ell(B+Z-D')=\ell(B-D')$.

Finally we see that this gives that $\rm{deg}(E) \geq \rm{deg}(Z)$.

\qed

We note that the above proof assumes that the code $C_{\Omega}(D,G)$ is 
nontrivial and thus has a codeword of nonzero weight.  If the code is 
trivial, the question of minimum distance is not very interesting.

\begin{rmrk}
Taking the special case $A = H$, $B = \lfloor H \rfloor$ and $Z = H - \lfloor 
H \rfloor$ in the AF bound gives the floor bound.  
Unlike the floor bound, the AF bound
can always be applied.  For an arbitrary divisor $G$ it is not always
the case that there exists a divisor $H$ such that $G$ can be written
as the sum of $H$ and $\lfloor H \rfloor$; however one can always let
$A=G$ and $B=Z=0$ and use the asymmetric bound trivially.  Beyond
this, there are many examples (see \S4) where the floor bound cannot
be applied, but the asymmetric bound gives an improvement over the
designed distance.  Additionally, there are examples where the floor
bound does apply, but the more general choice of $A$ and $B$ gives a
greater improvement.
\end{rmrk}

\begin{rmrk}
Setting $A = F + (\alpha + t)P$, $B = G + (\beta - t -1)P$ and $Z = (t + 1)P$ in the AF bound yields the K-P bound.  
So in the case of one point codes, the AF bound reduces to the K-P bound.  
Note that in its statement,
however, Theorem \ref{Prop3.10} makes no assumption about the support
of $F$ or $G$ containing only $P$.  Thus, the bound can be applied to
$n$-point codes, but only takes advantage of ``one-dimensional'' gaps,
while the asymmetric bound uses ``multi-dimensional'' gaps.  We also
note that in order to prove their bound, Kirfel and Pellikaan compare
their estimate to the F-R bound and show that their estimate is
always lower.  This shows that in the one point case, the asymmetric
bound will not improve on Feng-Rao.
\end{rmrk}

\section{Examples}
In this section we give several examples of how the asymmetric bound
improves over the floor bound and the K-P bound.

\begin{expl} \label{G41}
Consider the one point codes generated by the Suzuki curve, $y^8-y =
x^{10}-x^3$ over $\overline{\mathbb{F}}_8$.  This is a curve of genus
$g=14$.  If we let $G=41P_\infty$ and let $D$ be the sum of all the
other rational points on the curve, then we get that $C_\Omega(D,G)$
has designed distance 15.  From Table 2 in \cite{DuursmaChen} and the
equivalence $C_\Omega(D,aP_\infty)=C_\LL(D,(90-a)P_\infty)$ , we know
that the actual distance is 16.  The Weierstrass semi-group is
generated by 8, 10, 12, and 13, so we see that there is no way to
write $G$ as $H + \lfloor H \rfloor$, so that the floor bound does not
apply.  However, if we let $A=27P_\infty$, $B=14P_\infty$, and
$Z=P_\infty$ then Theorem \ref{asymbd} applies and we get that the
minimum distance is at least $15 + \rm{deg}(Z) = 16$, meeting the actual
distance.  Also, note that the K-P bound predicts this as
well since we are dealing with a one point code.
\end{expl}

\begin{expl} \label{G321}
This example illustrates the idea that the asymmetric bound can
exploit multiple ``multi-dimensional'' gaps to yield larger distance
estimates.  Consider the two point code, $C_\Omega(D,G)$, generated by the
same curve as in Example \ref{G41} where $G=32P_{0,0} + 1P_\infty$,
and $D$ is the sum of the remaining rational points.  The designed
distance of this code is seven.  By using knowledge of the Weierstrass
semigroup of the pair $(P_{0,0}, P_\infty)$ computed in
\cite{MatthewsS} we can apply all three bounds from the previous
section.  By letting $H=16P_{0,0} + 1P_\infty$, we can apply the floor
bound to get minimum distance at least eight.  Letting $S=1P_{0,0} +
0P_\infty$, $T=0P_{0,0} + 1P_\infty$, $\alpha=18$, $\beta =14$, and
$t=1$ shows that the K-P bound also gives minimum
distance at least eight.  However, using the asymmetric bound with
$A=14P_{0,0} + 1P_\infty$ and $B=18P_{0,0} + 0P_\infty$, we get
$Z=1P_{0,0} + 1P_\infty$ and a minimum distance of at least nine (See
Table 2).  
\end{expl}

\begin{expl}\label{G1517}
Consider the two point code, $C_\Omega(D,G)$, generated by the same
curve as in Example \ref{G41} where $G=15P_{0,0} + 17P_\infty$, and
$D$ is the sum of the remaining rational points.  If we let
$A=2P_{0,0} + 14P_\infty$ and $B=13P_{0,0} + 3P_\infty$, then
$Z=2P_{0,0} + 1P_\infty$ and the asymmetric bound gives distance at
least nine while $C_\Omega(D,G)$ has designed distance of six (See Table 2).  
The
floor bound, however, does not apply to this code or any linearly equivalent 
code.
\end{expl}

In the following we show the improvement over the designed distance 
achieved by the AF bound for two particular families of two-point codes; 
first the Hermitian codes over $\mathbb{F}_{16}$, then the Suzuki codes 
over $\mathbb{F}_8$.  \\ \\

\begin{tabular}{lcr}
\ &

\begin{tabular}{|c|ccccc|c|ccccc|c|ccccc|c|}
\hline
\ & 0 & 1 & 2 & 3 & 4 & \ & 0 & 1 & 2 & 3 & 4 & \  & 0 & 1 & 2 & 3 & 4 & \   \\
\hline
6 &  \ & \ & \ & \ & + & \ &  \ & \ & \ & \ & * & \ &  \ & \ & \ & \ & 2 & 6  \\
7 &  \ & \ & \ & + & + & \ &  \ & \ & \ & * & * & \ &  \ & \ & \ & 3 & 2 & 7  \\ 
8 &  \ & \ & + & + & + & \ &  \ & \ & 2 & * & * & \ &  \ & \ & 3 & 2 & 1 & 8  \\ 
9 &  \ & + & + & + & \ & \ &  \ & * & * & * & \ & \ &  \ & 2 & 1 & 1 & \ & 9  \\
10&  \ & + & + & + & \ & \ &  \ & 1 & 2 & * & \ & \ &  \ & 1 & 2 & 1 & \ & 10 \\ 
11&  + & + & + & + & + & \ &  1 & 2 & 3 & * & * & \ &  1 & 2 & 3 & 2 & 1 & 11 \\ 
12&  + & + & + & + & + & \ &  2 & 3 & * & * & * & \ &  2 & 3 & 2 & 1 & 1 & 12 \\ 
13&  + & + & + & \ & \ & \ &  * & * & * & \ & \ & \ &  1 & 2 & 1 & \ & \ & 13 \\ 
14&  \ & + & + & \ & \ & \ &  \ & * & * & \ & \ & \ &  \ & 1 & 1 & \ & \ & 14 \\ 
15&  \ & + & + & \ & \ & \ &  \ & * & * & \ & \ & \ &  \ & 1 & 1 & \ & \ & 15 \\
16&  + & + & + & + & + & \ &  * & 1 & * & * & * & \ &  1 & 1 & 1 & 1 & 1 & 16 \\
17&  + & + & \ & \ & \ & \ &  * & * & \ & \ & \ & \ &  1 & 1 & \ & \ & \ & 17 \\
18&  \ & + & \ & \ & \ & \ &  \ & 1 & \ & \ & \ & \ &  \ & 1 & \ & \ & \ & 18 \\
19&  \ & + & \ & \ & \ & \ &  \ & * & \ & \ & \ & \ &  \ & 1 & \ & \ & \ & 19 \\ 
20&  \ & + & \ & \ & \ & \ &  \ & 1 & \ & \ & \ & \ &  \ & 1 & \ & \ & \ & 20 \\
21&  + & \ & \ & \ & \ & \ &  1 & \ & \ & \ & \ & \ &  1 & \ & \ & \ & \ & 21 \\
\hline
\ & \ & \ & a & \ & \ & \ & \ & \ & b & \ & \ & \ & \  & \ & c & \ & \ & \ \\
\hline

\end{tabular} & \  \\ \\
\ & {Table 1 - Hermitian 2-point codes over $\mathbb{F}_{16}$ } & \  \\
\multicolumn{3}{l}{a. Codes where the actual distance is greater than the 
designed distance.}\\
\multicolumn{2}{l}{b. Improvement that the floor bound gives over the designed 
distance.}\\
\multicolumn{2}{l}{c. Improvement that the asymmetric bound gives over the designed 
distance.}\\
\end{tabular}

Let $\X$ be the Hermitian curve $y^4+y=x^5$ over $\F_{16}$, then
$g(\X) = 6$.  In Table~1a, a plus ``+'' designates those coordinates
$(x,y)$ such that the $C_\Omega(D,G)$ code with $G=xP_{0,0} + yP_\infty$,
$D$ the sum of the remaining rational points, and $\rm{deg}(G) \geq 10$ has
minimum distance strictly greater than the designed distance.  The coordinates
$(x,y)$ were calculated using Magma.  Any effective $G$
can be found in this table by use of the linear equivalence
$5P_{0,0}\sim 5P_\infty$.  In Tables 1b and 1c, the number $n$ in
position $(x,y)$ gives the improvement over the designed distance
predicted by the floor and asymmetric bounds, respectively.  An
asterisk~``*'' in Table~1b indicates that the floor bound either gave
no improvement, or could not be applied.  Note that Theorem
\ref{asymbd} gives an improvement for every code in Table 1a whose
actual distance is greater than its designed distance.  Also note that
more can be concluded from the information in Table~1b by considering
containment of codes.  For example, the code with $G=12P_{0,0} +
2P_\infty$ is contained in the code with $G=12P_{0,0} + 1P_\infty$,
thus its minimum distance is at least that of the larger code.

\begin{center}
\begin{tabular}{c}
\begin{tabular}{|c|ccccccccccccc|}

\hline
 \  & 0 & 1 & 2 & 3 & 4 & 5 & 6 & 7 & 8 & 9 & 10 & 11 & 12 \\
\hline 
14  & \ & \ & \ & \ & \ & \ & \ & \ & \ & \ & \  & \  & 2 \\
15  & \ & \ & \ & \ & \ & \ & \ & \ & \ & \ & \  & 3  & 3 \\
16  & \ & \ & \ & \ & \ & \ & \ & \ & \ & \ & 3  & 2  & 3 \\
17  & \ & \ & \ & \ & \ & \ & \ & \ & \ & 3 & 3  & 2  & 3 \\
18  & \ & \ & \ & \ & \ & \ & \ & \ & 3 & 3 & 3  & 2  & 2 \\
19  & \ & \ & \ & \ & \ & \ & \ & 4 & 3 & 3 & 3  & 2  & 2 \\
20  & \ & \ & \ & \ & \ & \ & 4 & 3 & 3 & 3 & 2  & 2  & 1 \\
21  & \ & \ & \ & \ & \ & 3 & 3 & 3 & 2 & 2 & 1  & 1  & \ \\
22  & \ & \ & \ & \ & 3 & 3 & 3 & 3 & 2 & 2 & 1  & 1  & 1 \\
23  & \ & \ & \ & 3 & 3 & 3 & 3 & 2 & 1 & 1 & \  & 1  & \ \\
24  & \ & \ & 3 & 2 & 2 & 2 & 2 & 2 & 1 & 1 & 1  & 1  & 1 \\
25  & \ & 2 & 3 & 3 & 3 & 2 & 2 & 1 & \ & 1 & \  & 1  & \ \\
26  & \ & 1 & 2 & 2 & 2 & 2 & 2 & 1 & \ & 1 & \  & 1  & \ \\
27  & 1 & 2 & 3 & 2 & 2 & 2 & 2 & 2 & 1 & 2 & 1  & 2  & 1 \\
28  & 2 & 3 & 4 & 3 & 3 & 3 & 2 & 3 & 2 & 2 & 2  & 1  & 1 \\
29  & 1 & 2 & 3 & 2 & 2 & 2 & 1 & 2 & 1 & 1 & 1  & \  & \ \\
30  & 1 & 2 & 3 & 2 & 2 & 2 & 2 & 2 & 1 & 1 & 1  & 1  & 1 \\
31  & 1 & 2 & 3 & 2 & 2 & 2 & 1 & 1 & 1 & \ & \  & \  & \ \\
32  & 1 & 2 & 2 & 1 & 2 & 1 & 2 & 1 & 1 & 1 & 1  & 1  & 1 \\
33  & 1 & 2 & 3 & 2 & 2 & 1 & 1 & \ & \ & \ & \  & \  & \ \\
34  & \ & 1 & 2 & 1 & 1 & 1 & 1 & \ & \ & \ & \  & \  & \ \\
35  & 1 & 2 & 2 & 1 & 1 & \ & 1 & \ & \ & \ & \  & \  & \ \\
36  & \ & 1 & 2 & 1 & 1 & \ & 1 & \ & \ & \ & \  & \  & \ \\
37  & 1 & 2 & 1 & \ & 1 & \ & 1 & \ & \ & \ & \  & \  & \ \\
38  & \ & 1 & 1 & \ & 1 & \ & 1 & \ & \ & \ & \  & \  & \ \\
39  & \ & 1 & 1 & \ & 1 & \ & 1 & \ & \ & \ & \  & \  & \ \\
40  & 1 & 1 & 1 & 1 & 1 & 1 & 1 & 1 & 1 & 1 & 1  & 1  & 1 \\
41  & 1 & 1 & 1 & 1 & 1 & 1 & \ & \ & \ & \ & \  & \  & \ \\
42  & \ & 1 & 1 & \ & 1 & \ & \ & \ & \ & \ & \  & \  & \ \\
43  & 1 & 1 & 1 & 1 & \ & \ & \ & \ & \ & \ & \  & \  & \ \\
44  & \ & 1 & 1 & \ & \ & \ & \ & \ & \ & \ & \  & \  & \ \\
45  & 1 & 1 & \ & \ & \ & \ & \ & \ & \ & \ & \  & \  & \ \\
46  & \ & 1 & \ & \ & \ & \ & \ & \ & \ & \ & \  & \  & \ \\
47  & \ & 1 & \ & \ & \ & \ & \ & \ & \ & \ & \  & \  & \ \\
48  & \ & 1 & \ & \ & \ & \ & \ & \ & \ & \ & \  & \  & \ \\
49  & \ & 1 & \ & \ & \ & \ & \ & \ & \ & \ & \  & \  & \ \\          
50  & \ & 1 & \ & \ & \ & \ & \ & \ & \ & \ & \  & \  & \ \\
51  & \ & 1 & \ & \ & \ & \ & \ & \ & \ & \ & \  & \  & \ \\
52  & \ & 1 & \ & \ & \ & \ & \ & \ & \ & \ & \  & \  & \ \\
53  & 1 & \ & \ & \ & \ & \ & \ & \ & \ & \ & \  & \  & \ \\
\hline

\end{tabular} \\
Table 2  \\ 
The asymmetric bound improvement for the Suzuki curve over $\F_8$.\\
\end{tabular}
\end{center}

Again consider the two-point Suzuki codes over $\mathbb{F}_8$ as in Examples 4.2 and 4.3.  
Table~2 gives the improvement over the designed distance of these codes using Theorem
\ref{asymbd}.  Table~2 is the analog for the Suzuki codes of 
Table 1c.  We list those divisors $G$ such that $\rm{deg} G \ge 2g - 2 = 26$ and such that
Theorem \ref{asymbd} predicts an improvement over the designed distance.

\begin{note}
We note that Tables 1 and 2 were computed using information about the
Weierstrass semigroup on two points that appeared in \cite{MatthewsW}
and \cite{MatthewsS} respectively.
\end{note}

We also note that the AF bound does not reach the actual distance in all cases.  Nor does it improve on the F-R bound for one-point codes.  In the general m-point case, both the floor bound and the AF bound produce improvements for the minimum distance yet they lack decoding algorithms to exploit those improvements.  

\section{Acknowledgments}
The authors would like to thank Professor Iwan Duursma for his advising and
tutelage through the course of this project.  Benjamin Lundell was
supported by an NSF VIGRE REU grant and Jason McCullough by an NSF VIGRE
graduate fellowship, both grant number DMS9983.

\nocite{*}
    
\bibliographystyle{plain} 

\bibliography{bib}

\end{document}